\documentclass[12pt,a4paper]{amsart}
\usepackage{amssymb}

\def\constr#1^#2{\mathrel{\mathop{\kern 0pt#1}\limits^{#2}}}
\def\build#1_#2{\mathrel{\mathop{\kern 0pt#1}\limits_{#2}}}
  
at 8pt

\begin{document}
\title[Grothendieck local duality and Cohomological Hasse principle]
{Grothendieck local duality and Cohomological Hasse principle for
$2$-dimensional complete local ring}

\author{Belgacem Draouil}

\address{D\'{e}partement de Math\'{e}matiques\\Facult\'{e} des
Sciences de Bizerte, 7021 Jarzouna}

\address{TUNISIA}

\email{Belgacem.Draouil@fsb.rnu.tn}

\date{November 2006}



\maketitle

\begin{abstract}
We prove a local duality for some schemes associated to a
$2$-dimensional complete local ring whose residue field is an
$n$-dimensional local field in the sense of Kato-Parshin. Our
results generalize the Saito works in the case $n=0$ and are applied
to study the Bloch-Ogus complex for such rings in various cases.

\end{abstract}

\section{ Introduction}

Let $A$ be a $2$-dimensional complete local ring with finite residue
field. The Bloch-Ogus complex associated to $A$ has been studied by
Saito in [12]. In this prospect, he calculated the homologies of
this complex and obtained (for any integer $n\geq1$) the following
exact sequence
{\footnotesize{
\begin{equation}
0\longrightarrow(\mathbb{Z}/n)^{r}\longrightarrow H^{3}(K,\mathbb{Z}
/n(2))\longrightarrow\underset{v\in P}{\oplus}
H^{2}(k(v),\mathbb{Z}/n(1))
\longrightarrow\mathbb{Z}/n\longrightarrow 0,
\tag{1.1}%
\end{equation}
}}
where $P$ denotes the set of height one prime ideals of $A$, $K~$is
the fractional field of $A$, $k(v)$ is the residue field at $v\in
P$, and $r=r(A)$ is an integer depending on the degeneracy of
$SpecA$ This result is based upon
the isomorphism ([12], lemma 5.4):%

\begin{equation}
H^{4}\left(  X,\mathbb{Z}/n\left(  2\right)  \right)
\simeq\mathbb{Z}/n,
\tag{1.2}%
\end{equation}
\ \ \ where $X=SpecA\backslash\left\{  x\right\}  $; and $x$ is the
unique maximal ideal of $A.$A decate later, Matsumi in [8]
generalised the result by Saito to $3$-dimensional complete regular
local ring of positive
characteristic. Indeed, he proved the exactness of the complex%

\begin{equation}%
\begin{array}
[c]{cc}%
0\longrightarrow H^{4}\left(  K,\mathbb{Z}/\ell\left(  3\right)
\right) \longrightarrow\underset{v\in\left(  SpecA\right)
_{2}}{\oplus}H^{3}\left( k\left(  v\right)  ,\mathbb{Z}/\ell\left(
2\right)  \right)  &\\
\longrightarrow\underset{v\in\left(
SpecA\right)  _{1}}{\oplus}H^{2}\left(
k\left(  v\right)  ,\mathbb{Z}/\ell\left(  1\right)  \right)
\longrightarrow\mathbb{Z}/\ell\longrightarrow0
\end{array}
\tag{1.3}%
\end{equation}
for all $\ell$ prime to $char(A)$, where $\left(  SpecA\right)
_{i}$ indicates the set of all points in $\ SpecA$ of dimension $i$.
Besides, if the ring $A$ is not regular, then the map \
\[
H^{4}\left(  K,\mathbb{Z}/\ell\left(  3\right)  \right)  \overset{\Psi_{K}%
}{\longrightarrow}\underset{v\in\left(  SpecA\right)  _{2}}{\oplus}%
H^{3}\left(  k\left(  v\right)  ,\mathbb{Z}/\ell\left(  2\right)
\right)
\]
\ is non-injective.

We proved in ([3], th3) that $Ker\Psi_{K}$ contains a sub-group of
type
$(\mathbb{Z}/\ell)^{r_{1}^{^{\prime}}(A)\text{ }}$, where $r_{1}^{^{\prime}%
}(A)$ is calculated as the $\mathbb{Z}-$rank of the graph of the
exceptional fiber of a resolution of $SpecA$. The main tools used in
this direction are
the isomorphism%

\begin{equation}
H^{6}\left(  X,\mathbb{Z}/\ell\left(  3\right)  \right)  \simeq\mathbb{Z}%
/\ell\tag{1.4}%
\end{equation}
and the perfect pairing%

\begin{equation}
H^{i}\left(  X,\mathbb{Z}/\ell\right)  \times H^{6-i}\left(  X,\mathbb{Z}%
/\ell(3)\right)  \longrightarrow H^{6}\left(
X,\mathbb{Z}/\ell\left(
3\right)  \right)  \simeq\mathbb{\ Z}/\ell\tag{1.5}%
\end{equation}
for all $i\geq1$ ([3], Section 0,D2).

In order to generalize the previous results, we consider a
$2$-dimensional complete local ring $A$ whose residue field $k$ is
$n$-dimensional local field in the sense of Kato-Parshin. To be more
precise, let $X=SpecA\backslash \left\{  x\right\}  $, where $x$ is
the unique maximal ideal of $A$. Then the main result of this paper
is the following

\textbf{Theorem (theorem 3.1)}

\textit{There exist an isomorphism }

\begin{equation}
H^{4+n}\left(  X,\mathbb{Z}/\ell\left(  2+n\right)  \right)  \simeq
\mathbb{Z}/\ell\tag{1.6}%
\end{equation}
\textit{and a perfect pairing }
{\footnotesize{
\begin{equation}
H^{1}\left(  X,\mathbb{Z}/\ell\right)  \times H^{3+n}\left(  X,\mathbb{Z}%
/\ell(2+n)\right)  \longrightarrow H^{4+n}\left(
X,\mathbb{Z}/\ell\left(
2+n\right)  \right)  \simeq\mathbb{\ Z}/\ell\tag{1.7}%
\end{equation}
}}
\textit{for all }$\ell$\textit{\ prime to }$char(A).$

We apply this result to calculate the homologies of the Bloch-Ogus
complex associated to $A.$ Indeed, let $\pi_{1}^{c.s}\left(
X\right)  $ be the quotient group of $\pi_{1}^{ab}\left(  X\right)
$ which classifies abelian c.s coverings of $X$ (see definition 4.1
below). We then prove the following

\textbf{Theorem (theorem 4.2)}

\textit{Let }$A$\textit{\ be a }$2$\textit{-dimensional complete
normal local
ring of positive characteristic whose residue field is }$n$%
\textit{-dimensional local field. Then the exact sequence }%

\begin{equation*}
0\longrightarrow\pi_{1}^{c.s}( X)/\ell \longrightarrow
H^{3+n}(K,\mathbb{Z}/\ell(2+n))\qquad\qquad
\end{equation*}
\begin{equation}
\qquad\qquad\qquad\longrightarrow\underset{v\in P}{\oplus}
H^{2+n}(k(v) ,\mathbb{Z}/\ell(1+n))
\longrightarrow\mathbb{Z}/\ell\longrightarrow0 \tag{1.8}%
\end{equation}

\bigskip\textit{holds.}

\textit{\ } $\,$Furthermore, if $A$ is assumed to be regular then,
as in the case of $n=0$ considered by Saito [12], we prove the
vanishing of the group $\pi_{1}^{c.s}\left(  X\right)  $ using the
recent paper [11] by Panin.

To prove these results, we rely heavily on the Grothendieck duality
theorem for strict local rings (section 3) as well as the purity
theorem of Fujiwara-Gabber, which we recall next.

In ([4], sentences just below Corollary7.1.7 ), Fujiwara confirmed
that the absolute cohomological purity in equicharacteristic is
true. In other words, we ge the following.

\textbf{Theorem of Fujiwara-Gabber}

\textit{Let }$T$\textit{\ be an equicharacteristic Notherian
excellent regular
scheme and }$Z$\textit{\ be a regular closed subscheme of codimension }%
$c$\textit{. Then for an arbitrary natural number }$\ell$\textit{\
prime to char(}$T$\textit{), the following canonical isomorphism }
\begin{equation}
H_{Z}^{i}\left(  T,\mathbb{Z}/\ell\left(  j\right)  \right)  \simeq
H^{i-2c}\left(  Z,\mathbb{Z}/\ell\left(  j-c\right)  \right)  \tag{1.9}%
\end{equation}

\bigskip\textit{holds.}

Finally, we complete the partial duality (1.7) in the case $n=1$ .
So, we obtain the following.

\textbf{Theorem (theorem 5.1)}

\textit{Let }$A$\textit{\ be a }$2$\textit{-dimensional normal
complete local
ring whose residue field is one-dimensional local field. Then, for very }%
$\ell$\textit{\ prime to }$char(A)$\textit{, the isomorphism }
\begin{equation}
H^{5}\left(  X,\mathbb{Z}/\ell\left(  3\right)  \right)  \simeq\mathbb{Z}%
/\ell\tag{1.10}%
\end{equation}
\textit{and the perfect pairing }
\begin{equation}
H^{i}\left(  X,\mathbb{Z}/\ell\right)  \times H^{5-i}\left(  X,\mathbb{Z}%
/\ell(3)\right)  \longrightarrow H^{5}\left(
X,\mathbb{Z}/\ell\left(
3\right)  \right)  \simeq\mathbb{\ Z}/\ell\tag{1.11}%
\end{equation}
\textit{hold for all }$i\in\{0,...,5\}$

\bigskip

Our paper is organised as follows. Section 2 devoted to some
notations. Section 3 contains the main theorem of this work
concerning the duality of the scheme $X=SpecA\backslash\left\{
x\right\}  $, where $A$ is a $2$-dimensional complete local ring
whose residue field is $n$-dimensional local field and $x$ is the
unique maximal ideal of $A.$In section 4, we study the Bloch-Ogus
complex associated to $A.$In section 5, we investigate the
particular case $n=1.$

\section{Notations}

For an abelian group $M$ and a positive integer $n$, we denote by
$M/n$ the cokernel of the map
$M\,\overset{n}{\longrightarrow}\,M.$For a scheme $Z$, and a sheaf
$\mathcal{F}$ over the \'{e}tale site of $Z$, $H^{i}\left(
Z,\mathcal{F}\right)  $ denotes the i-th \'{e}tale cohomology group.
For a positive integer $\ell$ invertible on $Z$, $\mathbb{\
Z}/\ell\left(  1\right) $ denotes the sheaf of $\ell$-th root of
unity and for an integer $i$, we denote $\mathbb{Z}/\ell\left(
i\right)  =\left(  \mathbb{\ Z}/\ell\left( 1\right)  \right)
^{\otimes i}$.

A local field $k$ is said to be $n$\textit{-dimensional local} if
there exists the following sequence of fields $k_{i}~\left(  1\leq
i\leq n\right)  $ such that

\noindent(i) each $k_{i}$ is a complete discrete valuation field
having $k_{i-1}$ as the residue field of the valuation ring
$O_{k_{i}}$ of $k_{i},$ and

\noindent(ii) $k_{0}$ is a finite field.

For such a field, and for $\ell$ prime to Char($k$), the well-known
isomorphism
\begin{equation}
H^{n+1}\left(  k,\mathbb{Z}/\ell\left(  n\right)  \right)  \simeq
\mathbb{Z}/\ell\tag{2.1}%
\end{equation}
and for each $i\in\{0,...,n+1\}$ a perfect duality%

\begin{equation}
H^{i}\left(  k,\mathbb{Z}/\ell\left(  j\right)  \right)  \times H^{n+1-i}%
\left(  k,\mathbb{Z}/\ell(n-j\right)  \longrightarrow H^{n+1}\left(
k,\mathbb{Z}/\ell\left(  n\right)  \right)  \simeq\mathbb{\ Z}/\ell\tag{2.2}%
\end{equation}

hold.

For a field $L$, $K_{i}\left(  L\right)  $ is the i-th Milnor group.
It coincides with the $i-$th Quillen group for $i\leq2.$ For $\ell$
prime to $char$ $L$, there is a Galois symbol
\begin{equation}
h_{\ell,L}^{i}\,\,\,\,K_{i}L/\ell\longrightarrow H^{i}(L,\mathbb{\ Z}%
/\ell\left(  i\right)  ) \tag{2.3}%
\end{equation}
which is an isomorphism for $i=0,1,2$ ($i=2$ is Merkur'jev-Suslin).
In this context, we recall the Kato conjecture ([7], Conjecture 1,
Section 1):

\textbf{Kato Conjecture}

\textit{For any field }$L$\textit{\ and any }$\ell$\textit{\ prime to }%
$charL$\textit{, the map }$h_{\ell,L}^{i}$\textit{\ is bijective.}

\section{Local duality}

We start this section by a description of the Grothendieck local
duality. Let $B$ denote a $d$-dimensional normal complete local ring
with maximal ideal $x^{\prime}$. By Cohen structure theorem ([10],
31.1), $B$ is a quotient of a regular local ring. Hence
\thinspace$SpecB$ admits a dualizing complex. Now, assume in the
first step that the residue field of $\,B$ is separably closed ($B$
is strictly local). Then, for
$X^{\prime}=SpecB\backslash\{x^{\prime}\}$ and for any $\ell$ prime
to char($A$), there is a Poincar\'{e} duality theory ([15],
Expos\'{e} I, Remarque 4.7.17). Namely, there is a trace isomorphism
\begin{equation}
H^{2d-1}\left(  X^{\prime},\mathbb{Z}/\ell\left(  d\right)  \right)
\overset{\backsim}{\longrightarrow}\mathbb{Z}/\ell\tag{3.1}%
\end{equation}
and a perfect pairing%

\begin{equation}
H^{i}\left(  X^{\prime},\mathbb{Z}/\ell\right)  \times
H^{2d-1-i}\left( X^{\prime},\mathbb{Z}/\ell(d\right)
\longrightarrow H^{2d-1}\left(
X^{\prime},\mathbb{Z}/\ell\left(  d\right)  \right)  \simeq\mathbb{\ Z}%
/\ell\tag{3.2}%
\end{equation}
\textit{for all }$i\in\{0,...,2d-1\}.$

Assume at this point that the residue field $k~$of $B$ is arbitrary.
Let $k_{s}$ be a separable closure of $k.$ The strict henselization
$B^{sh}$ of $B~$(with respect to the separably closed extension
$k_{s}$ of $k$) at the unique maximal ideal $x$ of $\ B$ is a
strictly local ring. It coincides with the integral closure of $B$
in the maximal unramified extension $L^{ur}$ of the fraction field
$L$ of $B$. Let $x^{\prime}$ be the maximal ideal of $B^{sh}$ and
let $X^{\prime}=SpecB^{sh}\backslash\{x^{\prime}\}$. So, the Galois
group of $\ X^{\prime}$ over $X$ is Gal($L^{ur}/L$) which is
isomorphic to Gal($k_{s}/k$). Then for any integer $j\geq0$, we get
the Hochschield-Serre
spectral sequence ([9], Remark 2.21)%

\begin{equation}
E_{2}^{p,q}=H^{p}(k,H^{q}(X^{\prime},\mathbb{Z}/\ell\left(  j\right)
)\Longrightarrow H^{p+q}(X,\mathbb{Z}/\ell\left(  j\right)  ) \tag{3.3}%
\end{equation}

Let $A$ denote a $2$-dimensional normal complete local ring whose
residue field is an $n$-dimensional local field. Let $x$ be the
unique maximal ideal of $A$. Then by normality $A$ admits at most
one singularity at $x$ in such a way that the scheme
$X=SpecA\backslash\left\{  x\right\}  $ becomes a regular scheme.

In what follows, we put

$K$ \ \ : the fractional field of $A$,

$k$ \ \ \ : the residue field of $K$

$P$ \ \ \ : the set of height one prime ideals of $A$.

For each $v\in P$ we denote by $K_{v}$ the completion of $K$ at $v$
and by $k(v)$ the residue field of $K_{v}$.

Let \ $X=SpecA\backslash\left\{  x\right\}  $ as above. Generalizing
(1.2), (1.4), and (1.5), we get the following.

\textbf{Theorem 3.1}

\textit{For all }$\ell$\textit{\ prime to }$char(A)$\textit{, the
isomorphism }
\begin{equation}
H^{4+n}\left(  X,\mathbb{Z}/\ell\left(  2+n\right)  \right)  \simeq
\mathbb{Z}/\ell\tag{3.4}%
\end{equation}
\textit{and the perfect pairing }
\begin{equation}
H^{1}\left(  X,\mathbb{Z}/\ell\right)  \times H^{3+n}\left(  X,\mathbb{Z}%
/\ell(2+n)\right)  \longrightarrow H^{4+n}\left(
X,\mathbb{Z}/\ell\left(
2+n\right)  \right)  \simeq\mathbb{\ Z}/\ell\tag{3.5}%
\end{equation}
\textit{occur. Furthermore, this duality is compatible with duality
(2.2) in the sense that the commutative diagram }
{\footnotesize{
\begin{equation}%
\begin{array}
[c]{cccccc}%
H^{1}\left(  X,\mathbb{Z}/\ell\right)  & \times & H^{3+n}\left(
X,\mathbb{Z}/\ell(2+n\right)  & \longrightarrow & H^{4+n}\left(
X,\mathbb{Z}/\ell\left(  2+n\right)  \right)  & \overset{\backsim
}{\longrightarrow}\mathbb{\ Z}/\ell\\
\downarrow i^{\ast} &  & \uparrow i_{\ast} &  & \uparrow i_{\ast} & \Vert\\
H^{1}\left(  k(v),\mathbb{Z}/\ell\right)  & \times & H^{n+1}\left(
k(v),\mathbb{Z}/\ell(n+1)\right)  & \longrightarrow & H^{n+2}\left(
k(v),\mathbb{Z}/\ell\left(  n+1\right)  \right)  & \overset{\backsim
}{\longrightarrow}\mathbb{\ Z}/\ell
\end{array}
\tag{3.6}%
\end{equation}
}}
\textit{holds, where }$\ i^{\ast}$\textit{\ is the map on }$\ H^{i}%
$\textit{\ induced from the map \ \ }$v\longrightarrow X$\textit{\ \
and }$i_{\ast}$\textit{\ \ is the Gysin map.}

\begin{proof}
The proof is slightly different from the proof of Theorem 1 in [2].
Let $k_{s}$ be a separable closure of $k.$We consider the strict
henselization $A^{sh}$ of $A$ (with respect to the separably closed
extension $k_{s}$ of $k$) at the unique maximal ideal $x$ of $\
A.~$Then, we denote $x^{\prime}$ the unique maximal ideal of
$A^{sh}$, $\ X^{\prime}=SpecA^{sh}\backslash \{x^{\prime}\}$ and we
use the spectral sequence (3.3). As $k$ is $n$-dimensional local
field, we have $H^{n+2}\left(  k,M\right)  =0\,\,$ for any torsion
module $M$ and as $X^{\prime}$ is of cohomological dimension $2d-1$
([14], the last paragraph of Introduction), we obtain
\[%
\begin{array}
[c]{ccc}%
H^{4+n}\left(  X,\mathbb{Z}/\ell(2+n)\right)  & \backsimeq & H^{n+1}%
(k,H^{3}(X^{\prime},\mathbb{Z}/\ell\left(  2+n\right)  )\\
& \backsimeq & H^{n+1}(k,\mathbb{Z}/\ell\left(  n\right)  )\text{ \
\ \ \ by
(3.1)}\\
& \backsimeq & \mathbb{Z}/\ell\text{
\ \ \ \ \ \ \ \ \ \ \ \ \ \ \ \ \ \ \ \ \ \ \ by (2.1)}%
\end{array}
\]

We prove now the duality (3.5). The filtration of the group
$H^{3+n}\left( X,\mathbb{Z}/\ell(2+n)\right)  $ is
\[
H^{3+n}\left(  X,\mathbb{Z}/\ell(2+n)\right)  =E_{n}^{3+n}\supseteq
E_{n+1}^{3+n}\supseteq0
\]
which leads to the exact sequence
\[
0\rightarrow E_{\infty}^{n+1,2}\longrightarrow H^{3+n}\left(  X,\mathbb{Z}%
/\ell(2+n)\right)  \longrightarrow E_{\infty}^{n,3}\longrightarrow0
\]

Since $E_{2}^{p,q}=0$ for all $p\geq n+2$ \ or $q\geq4\,$, we see
that
\[
E_{2}^{n,3}=E_{3}^{n,3}=...=E_{\infty}^{n,3}.
\]
\ \ \ \ \ \ \ \ \ \ \ \ \ The same argument yelds
\[
E_{3}^{n+1,2}=E_{4}^{n+1,2}=...=E_{\infty}^{n+1,2}%
\]
and $E_{3}^{n+1,2}=Co\ker d_{2}^{n-1,3}$ $\ $where $d_{2}^{n-1,3}$
is the map
\[
H^{n-1}(k,H^{3}(X^{\prime},\mathbb{Z}/\ell\left(  2+n\right)
)\longrightarrow H^{n+1}(k,H^{2}(X^{\prime},\mathbb{Z}/\ell\left(
2+n\right)  ).
\]
Hence, we obtain the exact sequence
\begin{equation}
0\rightarrow Co\ker d_{2}^{n-1,3}\longrightarrow H^{3+n}\left(  X,\mathbb{Z}%
/\ell(2+n)\right)  \longrightarrow H^{n}(k,H^{3}(X^{\prime},\mathbb{Z}%
/\ell\left(  2+n\right)  )\longrightarrow0 \tag{3.7}%
\end{equation}
Combining duality (2.2) for $k$ and duality (3.2), we deduce that
the group $H^{0}(k,H^{1}(X^{\prime},\mathbb{Z}/\ell))$ is dual to
the group $H^{n+1}(k,H^{2}(X^{\prime},\mathbb{Z}/\ell\left(
2+n\right)  )$ and the group
$H^{2}(k,H^{0}(X^{\prime},\mathbb{Z}/\ell))$ is dual to the group
$H^{n-1}(k,H^{3}(X^{\prime},\mathbb{Z}/\ell\left(  2+n\right)  )$ .
On the other hand, we have the commutative diagram
\begin{equation}%
\begin{array}
[c]{cccccc}%
H^{n-1}(k,H^{3}(X^{\prime},\mathbb{Z}/\ell\left(  2+n\right)  ) &
\times & H^{2}(k,H^{0}(X^{\prime},\mathbb{Z}/\ell)) &
\longrightarrow & H^{2}\left(
k,\mathbb{Z}/\ell\left(  1\right)  \right)  & \overset{\backsim}%
{\longrightarrow}\mathbb{Z}/\ell\\
\downarrow &  & \uparrow &  & \parallel & \Vert\\
H^{n+1}(k,H^{2}(X^{\prime},\mathbb{Z}/\ell\left(  2+n\right)  ) &
\times & H^{0}(k,H^{1}(X^{\prime},\mathbb{Z}/\ell)) &
\longrightarrow & H^{2}\left(
k,\mathbb{Z}/\ell\left(  1\right)  \right)  & \overset{\backsim}%
{\longrightarrow}\mathbb{Z}/\ell
\end{array}
\tag{3.8}%
\end{equation}
given by the cup products and the spectral sequence (3.3), using the
same argument as ( [1], diagram 46). We infer that $Co\ker
d_{2}^{n-1,3}$ is the dual of $Ker^{\prime}d_{2}^{0,1}$ where
$^{\prime}d_{2}^{0,1}$ is the boundary map for the spectral sequence
((3.4), j=0)
\begin{equation}
^{\prime}E_{2}^{p,q}=H^{p}(k,H^{q}(X^{\prime},\mathbb{Z}/\ell)\Longrightarrow
H^{p+q}(X,\mathbb{Z}/\ell) \tag{3.9}%
\end{equation}
Similarly, the group $H^{n}(k,H^{3}(X^{\prime},\mathbb{Z}/\ell\left(
2+n\right)  )$ is dual to the group $H^{1}(k,H^{0}(X^{\prime},\mathbb{Z}%
/\ell)).$The required duality is deduced from the following
commutative diagram
\[%
\begin{array}
[c]{ccccc}%
0\rightarrow Co\ker d_{2}^{n-1,3} & \longrightarrow & H^{3+n}\left(
X,\mathbb{Z}/\ell(2+n)\right)  & \longrightarrow &
H^{n}(k,H^{3}(X^{\prime
},\mathbb{Z}/\ell\left(  2+n\right)  )\longrightarrow0\\
\downarrow\wr &  & \downarrow\wr &  & \downarrow\wr\\
0\rightarrow(Ker^{\prime}d_{2}^{0,1})^{\vee} & \longrightarrow &
\left( H^{1}\left(  X,\mathbb{Z}/\ell\right)  \right)  ^{\vee} &
\longrightarrow &
(H^{1}(k,H^{0}(X^{\prime},\mathbb{Z}/\ell)))^{\vee}\longrightarrow0
\end{array}
\]
where the upper exact sequence is (3.7) and the bottom exact
sequence is the dual of the well-known exact sequence
\[
0\rightarrow^{\prime}E_{2}^{1,0}\longrightarrow H^{1}\left(  X,\mathbb{Z}%
/\ell\right)  \longrightarrow
Ker^{\prime}d_{2}^{0,1}\longrightarrow0
\]
deduced from the spectral sequence (3.9) and where $\left(  M\right)
^{\vee}$ denotes the dual $Hom(M,\mathbb{Z}/\ell)$ for any
$\mathbb{Z}/\ell-$module $M.$

Finally, to obtain the last part of the theorem, we remark that the
commutativity of the diagram (3.6) is obtained by via a same
argument (projection formula ([9], VI 6.5) and compatibility of
traces ([9], VI 11.1)) as [1] to establish the commutative diagram
in the proof of assertion ii) at page 791.
\end{proof}

\bigskip

\textbf{Corollary 3.2}

\textit{With the same notations as above, the following commutative
diagram }
\begin{equation}%
\begin{array}
[c]{ccc}%
H^{n+1}\left(  k(v),\mathbb{Z}/\ell(n+1)\right)  & \overset{i_{\ast}%
}{\longrightarrow} & H^{3+n}\left(  X,\mathbb{Z}/\ell(2+n\right) \\
\downarrow &  & \downarrow\\
(H^{1}\left(  k(v),\mathbb{Z}/\ell\right)  )^{\vee} &
\overset{(i^{\ast })^{\wedge}}{\longrightarrow} & (H^{1}\left(
X,\mathbb{Z}/\ell\right)
)^{\vee}%
\end{array}
\tag{3.10}%
\end{equation}

\textit{holds.}

\begin{proof}
\bigskip This is a consequence of diagram (3.6).
\end{proof}

\bigskip

The duality (3.5) will be completed (section 5) \ to a general
pairing by replacing $H^{1}\left(  X,\mathbb{Z}/\ell\right)  $ \ by
$H^{i}\left( X,\mathbb{Z}/\ell\right)  $ ; for $0\leq i\leq5$ in the
case $n=1.$

\section{The Bloch-Ogus complex}

In this section, we investigate the study of the Bloch-Ogus complex
associated to the ring $A~$considered previously. So, let $A$ be a
$2$-dimensional normal complete local ring whose residue field is an
$n$-dimensional local field. Next, we define a group which appears
in the homologies of the associated Bloch-Ogus complex of $A$.

\textbf{Definition 4.1}

\textit{Let }$Z$\textit{\ be a Noetherian scheme. A finite etale
covering }$\ f:$\textit{\ }$W\rightarrow Z$\textit{\ is called a c.s
covering if for any closed point }$z$\textit{\ of }$Z$\textit{,
}$z\times_{Z}W$\textit{\ is isomorphic to a finite scheme-theoretic
sum of copies of }$z$\textit{. We
denote }$\pi_{1}^{c.s}\left(  Z\right)  $\textit{\ the quotient group of }%
$\pi_{1}^{ab}\left(  Z\right)  $\textit{\ which classifies abelian
c.s coverings of }$Z.$

\bigskip

As above, let \ $X=SpecA\backslash\left\{  x\right\}  $. The group
$\pi _{1}^{c.s}\left(  X\right)  /\ell$ is the dual of the kernel of
the map
\begin{equation}
H^{1}\left(  X,\mathbb{Z}/\ell\right)  \longrightarrow%
{\displaystyle\prod\limits_{v\in P}}
H^{1}\left(  k\left(  v\right)  ,\mathbb{Z}/\ell\right)  \tag{4.1}%
\end{equation}
( [12], section 2, definition and sentence just below). Now, we are
able to calculate the homologies of the Bloch-Ogus complex
associated to the ring $\ A.$

\textbf{Theorem 4.2}

\textit{\ For all }$\ell$\textit{\ prime to the characteristic of }%
$A$\textit{, the following sequence is exact.}
\begin{equation}
0\longrightarrow\pi_{1}^{c.s}\left(  X\right)  /\ell\longrightarrow
H^{n+3}\left(  K,\mathbb{Z}/\ell\left(  n+2\right)  \right)
\longrightarrow
\underset{v\in P}{\oplus}H^{n+2}\left(  k\left(  v\right)  ,\mathbb{Z}%
/\ell\left(  n+1\right)  \right)  \longrightarrow\mathbb{\ Z}/\ell
\longrightarrow0 \tag{4.2}%
\end{equation}

\begin{proof}
Consider the localisation sequence on $X=SpecA\backslash\left\{  x\right\}  $%
\[
...\rightarrow H^{i}\left(  X,\mathbb{Z}/\ell(n+2)\right)
\longrightarrow H^{i}\left(  K,\mathbb{Z}/\ell(n+2)\right)
\longrightarrow\underset{v\in P}{\oplus}H_{v}^{i+1}\left(
X,\mathbb{Z}/\ell(n+2)\right)  \rightarrow...
\]
Firstly, for any $v\in P$ , we have the isomorphisms%

\[
H_{v}^{i}\left(  X,\mathbb{Z}/\ell\left(  2+n\right)  \right)
\simeq H_{v}^{i}\left(  SpecA_{v},\mathbb{Z}/\ell\left(  2+n\right)
\right)  \text{
\ \ \ \ }%
\]
by excision. Secondly, we can apply the purity theorem of
Fujiwara-Gabber (Introduction) for $Z=v,T=SpecA_{v}$ and we find the
isomorphisms
\[
H_{v}^{3+n}\left(  SpecA_{v},\mathbb{Z}/\ell\left(  2+n\right)
\right) \simeq H^{1+n}\left(  k\left(  v\right)
,\mathbb{Z}/\ell\left(  1+n\right) \right)  \
\]
and
\[
H_{v}^{4+n}\left(  SpecA_{v},\mathbb{Z}/\ell\left(  2+n\right)
\right) \simeq H^{2+n}\left(  k\left(  v\right)
,\mathbb{Z}/\ell\left(  1+n\right) \right)  \
\]
which lead to the isomorphisms%

\[
H_{v}^{3+n}\left(  X,\mathbb{Z}/\ell\left(  2+n\right)  \right)
\simeq H^{1+n}\left(  k\left(  v\right)  ,\mathbb{Z}/\ell\left(
1+n\right)  \right)
\]

and
\[
H_{v}^{4+n}\left(  X,\mathbb{Z}/\ell\left(  2+n\right)  \right)
\simeq H^{2+n}\left(  k\left(  v\right)  ,\mathbb{Z}/\ell\left(
1+n\right)  \right) .\
\]
Hence we derive the exact sequence

$\underset{v\in P}{\oplus}H^{1+n}\left(  k\left(  v\right)  ,\mathbb{Z}%
/\ell\left(  1+n\right)  \right)
\overset{g}{\longrightarrow}H^{3+n}\left(
X,\mathbb{Z}/\ell\left(  2+n\right)  \right)  \longrightarrow H^{3+n}%
(K,\mathbb{Z}/\ell(2+n))$

$\ \ \ \ \ \ \ \ \ \ \ \ \ \ \ \longrightarrow\underset{v\in P}{\oplus}%
H^{2+n}\left(  k\left(  v\right)  ,\mathbb{Z}/\ell\left(  1+n\right)
\right) \longrightarrow H^{4+n}\left(  X,\mathbb{Z}/\ell\left(
2+n\right)  \right) \longrightarrow0$

The last zero on the right is a consequence of the vanishing of the
group

$H^{4+n}(K,\mathbb{Z}/\ell(2+n))$. Indeed, $A$ is finite over
$O_{L}[[T]]$ for some complete discrete valuation field $L$ having
the same residue field with $A$ [10, \S 31]. By Serre [13, chapI,
Prop 14], $cd_{\ell}(A)\leq cd_{\ell }(O_{L}[[T]])$ and by Gabber
[5], $cd_{\ell}(O_{L}[[T]])=n+3\,\,\,$using the fact that
$cd_{\ell}(k)=n+1.$

Now, by the right square of the diagram (3.6), the Gysin map
\[
\underset{v\in P}{\oplus}H^{2+n}\left(  k\left(  v\right)  ,\mathbb{Z}%
/\ell\left(  1+n\right)  \right)  \longrightarrow H^{4+n}\left(
X,\mathbb{Z}/\ell\left(  2+n\right)  \right)
\]
can be replaced by the map $\underset{v\in P}{\oplus}H^{2+n}\left(
k\left( v\right)  ,\mathbb{Z}/\ell\left(  1+n\right)  \right)
\longrightarrow \mathbb{Z}/\ell$ after composing with the trace
isomorphism $H^{4+n}\left( X,\mathbb{Z}/\ell\left(  2+n\right)
\right)  \simeq\mathbb{Z}/\ell$ \ (3.4). So, we obtain the exact
sequence
\[
0\longrightarrow Co\ker g\longrightarrow H^{3+n}(K,\mathbb{Z}/\ell
(2+n))\longrightarrow\underset{v\in P}{\oplus}H^{2+n}\left(  k\left(
v\right)  ,\mathbb{Z}/\ell\left(  1+n\right)  \right)
\longrightarrow \mathbb{Z}/\ell\longrightarrow0
\]

Finally, in view of the commutative diagram (3.10), we deduce that
$Co\ker g$ equals to the group $\pi_{1}^{c.s}\left(  X\right)
/\ell$ taking in account (4.1).$\ \ \ \ \ \ \ \ \ \ \ \ \ $
\end{proof}

\bigskip

Next, we assume further that $A$ is regular. We will prove that the
group $\pi_{1}^{c.s}\left(  X\right)  /\ell$ \ vanishes.

\bigskip

\textbf{Theorem 4.3}

\textit{Let }$A=F_{p}((t_{1}))((t_{2}))...((t_{n}))[[X,Y]]$\textit{\
of fraction field }$K$\textit{\ and assum Kato conjecture (section
2), then the following Hasse principle complex for }$K$\textit{\ }
\begin{equation}
\,0\longrightarrow H^{3+n}(K,\mathbb{Z}/\ell(2+n))\longrightarrow
\underset{v\in P}{\oplus}H^{2+n}\left(  k\left(  v\right)  ,\mathbb{Z}%
/\ell\left(  1+n\right)  \right)  \longrightarrow\mathbb{Z}/\ell
\longrightarrow0 \tag{4.3}%
\end{equation}
\textit{is exact.}

\begin{proof}
Keeping in mind (4.2), \noindent it remains to prove the injectivity
of the map
\[
\Psi_{K}:\,\,\,\,\,\,H^{3+n}\left(  K,\mathbb{Z}/\ell\left(
2+n\right) \right)  \longrightarrow\underset{v\in
P}{\oplus}H^{2+n}\left(  k\left( v\right)  ,\mathbb{Z}/\ell\left(
1+n\right)  \right)  .
\]

Let $\ q$ be an integer and consider the sheaf
$\mathcal{H}^{q}\left( \mathbb{Z}/\ell\left(  n+2\right)  \right)  $
on $SpecA$, the Zariskien sheaf
associated to the presheaf $U\longrightarrow H^{q}\left(  U,\mathbb{Z}%
/\ell\left(  n+2\right)  \right)  $. As a consequence of a recent
work of Panin [11], we conclude that the cohomology of this sheaf is
calculated as the homology of the Bloch-Ogus complex, that is :
\[
H^{q}\left(  K,\mathbb{Z}/\ell\left(  n+2\right)  \right)
\longrightarrow
\underset{v\in P}{\oplus}H^{q-1}\left(  k\left(  v\right)  ,\mathbb{Z}%
/\ell\left(  n+1\right)  \right)  \longrightarrow H^{q-1}\left(
k\left( x\right)  ,\mathbb{Z}/\ell\left(  n\right)  \right)  .
\]
So the group $Ker\Psi_{K}$ is identified with the group $H^{0}\left(
(SpecA)_{Zar},\mathcal{H}^{n+3}(\mathbb{Z}/\ell\left(  n+2\right)
\right) ).$

On the other hand, the Bloch-Ogus spectral sequence
\[
H^{p}\left(  (SpecA)_{Zar},\mathcal{H}^{q}(\mathbb{Z}/\ell\left(
n+2\right) \right)  )\Rightarrow H^{p+q}(SpecA,\mathbb{Z}/\ell\left(
n+2\right)  )
\]
gives the exact sequence\newline$H^{n+3}(SpecA,\mathbb{Z}/\ell\left(
n+2\right)  )\longrightarrow H^{0}\left(  (SpecA)_{Zar},\mathcal{H}%
^{n+3}(\mathbb{Z}/\ell\left(  n+2\right)  \right)  )$

$\ \ \ \ \ \ \ \ \ \ \ \ \ \ \ \ \ \longrightarrow H^{2}\left(  (SpecA)_{Zar}%
,\mathcal{H}^{n+2}(\mathbb{Z}/\ell\left(  n+2\right)  \right)
)\rightarrow H^{n+4}(SpecA,\mathbb{Z}/\ell\left(  n+2\right)  )$

Sinces the ring $A$ is henselian, we obtain the isomorphism
\[
H^{i}\left(  SpecA,\mathbb{Z}/\ell(n+2)\right)  \simeq H^{i}\left(
Speck,\mathbb{Z}/\ell(n+2)\right)  ,i\geq0.
\]
But the groups $H^{n+3}\left(  Speck,\mathbb{Z}/\ell(3)\right)  $
and $H^{n+4}\left(  Speck,\mathbb{Z}/\ell(3)\right)  $ vanish
because the cohomological dimension of $k$ is $n+1.$ Thus we get the
isomorphism
\[
H^{0}\left(  (SpecA)_{Zar},\mathcal{H}^{n+3}(\mathbb{Z}/\ell\left(
n+2\right)  \right)  )\overset{\backsim}{\longrightarrow}H^{2}\left(
(SpecA)_{Zar},\mathcal{H}^{n+2}(\mathbb{Z}/\ell\left(  n+2\right)
\right)  )
\]
which means that \ $Ker\Psi_{K}$ is isomorphic to the Cokernel of
the map
\[
\underset{v\in P}{\oplus}H^{n+1}\left(  k\left(  v\right)  ,\mathbb{Z}%
/\ell\left(  n+1\right)  \right)  \longrightarrow H^{n}\left(
k\left( x\right)  ,\mathbb{Z}/\ell\left(  n\right)  \right)
\]
So, we must prove the surjectivity of this last map. Indeed, the
Gersten-Quillen complex ([11], Theorem A)
\[
K_{n+2}\left(  A\right)  /\ell\longrightarrow K_{n+2}\left(
K\right) /\ell\longrightarrow\underset{v\in P}{\oplus}K_{n+1}k\left(
v\right) /\ell\longrightarrow K_{n}k\left(  x\right)
/\ell\longrightarrow0
\]
is exact. On the other hand, we have the following commutative
diagram
\[%
\begin{array}
[c]{ccccc}%
K_{n+2}K/\ell & \longrightarrow & \underset{v\in
P}{\oplus}K_{n+1}k\left( v\right)  /\ell & \longrightarrow &
K_{n}k\left(  x\right)  /\ell
\longrightarrow0\\
\downarrow &  & \downarrow &  & \downarrow\wr\\
H^{n+2}\left(  K,\mathbb{Z}/\ell\left(  n+2\right)  \right)  &
\longrightarrow
& \underset{v\in P}{\oplus}H^{n+1}\left(  k\left(  v\right)  ,\mathbb{Z}%
/\ell\left(  n+1\right)  \right)  & \longrightarrow & H^{n}\left(
k\left( x\right)  ,\mathbb{Z}/\ell\left(  n\right)  \right)
\end{array}
\]
where the right vertical isomorphism comes from Kato conjecture.
This yelds that the map
\[
\underset{v\in P}{\oplus}H^{n+1}\left(  k\left(  v\right)  ,\mathbb{Z}%
/\ell\left(  n+1\right)  \right)  \longrightarrow H^{n}\left(
k\left( x\right)  ,\mathbb{Z}/\ell\left(  n\right)  \right)
\]
is surjective and we are done.
\end{proof}

\bigskip

\textbf{Remark 4.4}

1) The case $n=0$, implies the following exact sequence
\[
\,0\longrightarrow
H^{3}(K,\mathbb{Z}/\ell(2))\longrightarrow\underset{v\in
P}{\oplus}H^{2}\left(  k\left(  v\right)  ,\mathbb{Z}/\ell\left(
1\right) \right)  \longrightarrow\mathbb{Z}/\ell\longrightarrow0
\]
already obtained by Saito [12].

2) The case $n=1$ leads to the following exact sequence
\[
\,0\longrightarrow
H^{4}(K,\mathbb{Z}/\ell(3))\longrightarrow\underset{v\in
P}{\oplus}H^{3}\left(  k\left(  v\right)  ,\mathbb{Z}/\ell\left(
2\right) \right)  \longrightarrow\mathbb{Z}/\ell\longrightarrow0
\]
which is considered in [1].

3) The case $n=2.$ Let $A=\mathbb{F}_{p}((t))((u))[[X,Y]]$ of
fraction field $K$. Then the following Hasse principle complex for
$K$
\[
\,0\longrightarrow
H^{5}(K,\mathbb{Z}/\ell(4))\longrightarrow\underset{v\in
P}{\oplus}H^{4}\left(  k\left(  v\right)  ,\mathbb{Z}/\ell\left(
3\right) \right)  \longrightarrow\mathbb{Z}/\ell\longrightarrow0
\]
is exact.

4) The case $n=0$ has been used by Saito to study the class field
theory of curves over one dimensional local field. Recently, Yoshida
[16] provided an alternative approch which includes the equal
characteristic case. In a forthcoming paper I use the case $n=1$ to
investigate the study of class field theory of curves over $2-$
dimensional local field.

\section{The case \ n=1}

Let $A$ denote a $2$-dimensional complete normal local ring of
positive characteristic whose residue field is one-dimensional local
field. The aim of this section is to complete the duality (3.5) for
$i\geq1.$ We prove the following.

\textbf{Theorem 5.1 }

\textit{For all }$\ell$\textit{\ prime to }$char(A)$\textit{, the
isomorphism }
\begin{equation}
H^{5}\left(  X,\mathbb{Z}/\ell\left(  3\right)  \right)  \simeq\mathbb{Z}%
/\ell\tag{5.1}%
\end{equation}
\textit{and the perfect pairing }
\begin{equation}
H^{i}\left(  X,\mathbb{Z}/\ell\right)  \times H^{5-i}\left(  X,\mathbb{Z}%
/\ell(3)\right)  \longrightarrow H^{5}\left(
X,\mathbb{Z}/\ell\left(
3\right)  \right)  \simeq\mathbb{\ Z}/\ell\tag{5.2}%
\end{equation}
\textit{hold for all }$i\in\{0,...,5\}.$

\begin{proof}
The first isomorphism is given by (3.4). Next. We proceed to the
second part of the theorem. As in the proof of Theorem 3.1, we
consider the strict henselization $A^{sh}$ of $A$ (with respect to
the separably closed extension $k_{s}$ of $k$) at $x$. If
$x^{\prime}$ is the maximal ideal of $A^{sh}$, we recall that we
denote $X^{\prime}=SpecA^{sh}\backslash\{x^{\prime}\}$ and we
consider the spectral sequence \ (3.3). The filtration of the group
$H^{i}\left(  X,\mathbb{Z}/\ell(3)\right)  $ is
\[
H^{i}\left(  X,\mathbb{Z}/\ell(3)\right)  =E_{0}^{i}\supseteq E_{1}%
^{i}\supseteq E_{2}^{i}\supseteq0
\]
where the quotients are given by
\[
E_{0}^{i}/E_{1}^{i}\simeq E_{\infty}^{0,i}\simeq\ker d_{2}^{0,i}%
\]%
\[
E_{1}^{i}/E_{2}^{i}\simeq E_{\infty}^{1,i-1}\simeq E_{2}^{1,i-1},\text{ and}%
\]%
\[
\text{ \ }E_{2}^{i}\simeq E_{\infty}^{2,i-2}\simeq Co\ker d_{2}^{0,i-1}%
\]
$.$

The same computation is true for the group $H^{5-i}\left(  X,\mathbb{Z}%
/\ell\right)  $ and the filtration
\[
H^{5-i}\left(  X,\mathbb{Z}/\ell\right)
=^{\prime}E_{0}^{5-i}\supseteq
^{\prime}E_{1}^{5-i}\supseteq^{\prime}E_{2}^{5-i}\supseteq0
\]
by considering the spectral sequence ((3.4), j=0).

Now, combining duality (2.2) and duality (3.2) we observe that the
group $E_{2}^{0,j}$ is dual to the group $^{\prime}E_{2}^{2,3-j}$
and the group $E_{2}^{1,j}$ is dual to the group
$^{\prime}E_{2}^{1,3-j}$ for all \ $0\leq j\leq3.$ On the other
hand, we have the commutative diagram
\[%
\begin{array}
[c]{cccccc}%
H^{0}(k,H^{i}(X^{\prime},\mathbb{Z}/\ell\left(  3\right)  ) & \times
& H^{2}(k,H^{3-i}(X^{\prime},\mathbb{Z}/\ell)) & \longrightarrow &
H^{2}\left(
k,\mathbb{Z}/\ell\left(  1\right)  \right)  & \overset{\backsim}%
{\longrightarrow}\mathbb{Z}/\ell\\
\downarrow &  & \uparrow &  & \parallel & \Vert\\
H^{2}(k,H^{i-1}(X^{\prime},\mathbb{Z}/\ell\left(  3\right)  ) &
\times & H^{0}(k,H^{4-i}(X^{\prime},\mathbb{Z}/\ell)) &
\longrightarrow & H^{2}\left(
k,\mathbb{Z}/\ell\left(  1\right)  \right)  & \overset{\backsim}%
{\longrightarrow}\mathbb{Z}/\ell
\end{array}
\]
given by the cup products and the spectral sequence (3.4), using the
same argument as ( [1],diagram 46). We infer that $Co\ker
d_{2}^{0,i-1}$ is the dual of $Ker^{\prime}d_{2}^{0,5-i}$ and
$Kerd_{2}^{0,i}$ is the dual of $Co\ker^{\prime}d_{2}^{0,2d-i}$where
$^{\prime}d_{2}^{p,q}$ is the boundary map for the spectral sequence
((3.4),j=0)
\[
^{\prime}E_{2}^{p,q}=H^{p}(k,H^{q}(X^{\prime},\mathbb{Z}/\ell)\Longrightarrow
H^{p+q}(X,\mathbb{Z}/\ell).
\]

This is illustrated by the following diagram:%

\begin{equation}%
\begin{array}
[c]{ccccc}%
H^{i}X,\mathbb{Z}/\ell(3))=E_{0}^{i} &  &  &  & 0\\
\mid & \left\}  Kerd_{2}^{0,i}\right.  & \longleftrightarrow &
\left.
Co\ker^{\prime}d_{2}^{0,24-i}\right\{  & \mid\\
E_{1}^{i} &  &  &  & ^{\prime}E_{2}^{5-i}\\
\mid & \left\}  E_{2}^{1,i-1}\right.  & \longleftrightarrow & \left.
^{\prime}E_{2}^{1,24-i}\right\{  & \mid\\
E_{2}^{i} &  &  &  & ^{\prime}E_{1}^{5-i}\\%
\begin{array}
[c]{c}%
\mid\\
0
\end{array}
& \left\}  Co\ker d_{2}^{0,i-1}\right.  & \longleftrightarrow &
\left. Ker^{\prime}d_{2}^{0,5-i}\right\{  &
\begin{array}
[c]{c}%
\mid\\
^{\prime}E_{0}^{5-i}=H^{5-i}\left(  X,\mathbb{Z}/\ell\right)
\end{array}
\end{array}
\tag{5.3}%
\end{equation}

where each pair of groups which are combined by
$\longleftrightarrow$ consists of a group and its dual group.

We begin by calculate the dual group of $E_{1}^{i}$, using the
following
commutative diagram%

\[%
\begin{array}
[c]{ccccccc}%
0\rightarrow & E_{2}^{i} & \rightarrow & E_{1}^{i} & \rightarrow & E_{1}%
^{i}/E_{2}^{i} & \rightarrow0\\
& \downarrow\wr &  & \downarrow &  & \downarrow\wr & \\
0\rightarrow &
(^{\prime}E_{0}^{2d+1-i}/^{\prime}E_{1}^{2d+1-i})^{\vee} &
\rightarrow &
(^{\prime}E_{0}^{2d+1-i}/^{\prime}E_{2}^{2d+1-i})^{\vee} &
\rightarrow &
(^{\prime}E_{1}^{2d+1-i}/^{\prime}E_{2}^{2d+1-i})^{\vee} &
\rightarrow0
\end{array}
\]

where $\left(  M\right)  ^{\vee}$ denotes the dual
$Hom(M,\mathbb{Z}/\ell)$ for any $\mathbb{Z}/\ell-$module $M$ and
where the left and right vertical isomorphisms are explained by the
previous diagram. This yelds that
\begin{equation}
E_{1}^{i}\simeq(^{\prime}E_{0}^{5-i}/^{\prime}E_{2}^{5-i})^{\vee} \tag{5.4}%
\end{equation}

Finally, the required duality between
$H^{i}X,\mathbb{Z}/\ell(d+1))=E_{0}^{i}$ \ and \ \
$^{\prime}E_{0}^{5-i}=H^{5-i}\left(  X,\mathbb{Z}/\ell\right)  $
follows from the following commutative diagram
\[%
\begin{array}
[c]{ccccccc}%
0\rightarrow & E_{1}^{i} & \rightarrow & E_{0}^{i} & \rightarrow & E_{0}%
^{i}/E_{1}^{i} & \rightarrow0\\
& \downarrow\wr &  & \downarrow\wr &  & \downarrow\wr & \\
0\rightarrow & (^{\prime}E_{0}^{5-i}/^{\prime}E_{2}^{5-i})^{\vee} &
\rightarrow & (^{\prime}E_{0}^{5-i})^{\vee} & \rightarrow & (^{\prime}%
E_{2}^{5-i})^{\vee} & \rightarrow0
\end{array}
\]

where the right vertical isomorphism is given by (5.3) and the left
vertical isomorphism is the isomorphism (5.4).
\end{proof}

$\bigskip$


\begin{thebibliography}{99}                                                                                               %


\bibitem {1} Colliot-Th\'{e}l\`{e}ne,J.L.,Sansuc,J.J.,Soul\'{e},C. \ Torsion
dans le groupe de Chow de codimension deux. \ Duke Math. Journal
vol. 50 No.3 pp763-801 (1983)

\bibitem {2} Draouil,B. Cohomological Hasse principle for the ring
$\ \mathbb{F}_{p}((t))[[X,Y]]$, Bull. Belg. Math. Soc. Simon Stevin
11, no. 2 (2004), 181--190

\bibitem {3} Draouil, B.,Douai, J.C. \thinspace\ Sur l'arithm\'{e}tique des
anneaux locaux de dimension 2 et 3, Journal of Algebra 213 (1999),
499-512.

\bibitem {4} Fujiwara, K. $\,\,\,$Theory of Tubular Neighborhood in Etale
Topology Duke Math.J.80 (1995), 15-57.

\bibitem {5} Gabber,O. \ Lectuure at IHES,on March 1981.

\bibitem {6} Kato, K A Hasse principle for two-dimensional global fields,
J.reine angew.Math.366 (1986), 143-183.

\bibitem {7} Kato, K. \thinspace\thinspace\thinspace\thinspace A
genralisation of local class field theory by using K-theory II,
J.Fac.Sci.Univ. Tokyo, 27 (1980), 603-683.

\bibitem {8} Matsumi, K.\thinspace\thinspace\ Thesis, Arithmetic of
three-dimensional complete regular local rings of positive
characteristics
T\^{o}%
${{}^2}$%
hoku University, Japan 1999.

\bibitem {9} Milne;J.S. \ Etale Cohomology, \ Princeton University Press,
Princeton 1980.

\bibitem {10} Nagata,M. \ Local rings,Tracts in Mathematics Number 13,
Intersciences Publishers. New York 1962\ \ \ \

\bibitem {11} Panin,I. \ The equi-characteristic case of the \ Gersten
conjecture, preprint (2000),available on the K-theory server.

\bibitem {12} Saito, S. $\,\,\,$Class field Theory for two-dimensional
local rings Galois groups and their representations,
Kinokuniya-North Holland Amsterdam, vol 12 (1987), 343-373

\bibitem {13} Serre,J.P. \ \ Cohomologie Galoisienne ,L.N.M 5
\ Berlin-Heidelberg-New York 1965.

\bibitem {14} SGA 4 \ \thinspace\thinspace\thinspace Th\'{e}orie des Topos
et Cohomologie \'{e}tale des Sch\'{e}mas. Lecture Notes in Math.
vol.305, Springer-Verlag \ Berlin.Heidelberg.New York.

\bibitem {15} SGA 5 \ \thinspace\thinspace\thinspace Cohomologie l-adique
et Fonctions L . \ Lecture Notes in Math. vol.589, Springer-Verlag \
Berlin.Heidelberg.New York

\bibitem {16} Yoshida. Finitness theorems in the class field theory of
varieties over local fields, Journal of Number Theory 101 (2003)
138-150
\end{thebibliography}
\end{document}